\documentclass[12pt,reqno]{article}
 \usepackage{amssymb}
 \usepackage{hyperref}
 \usepackage{amsmath}
 \usepackage{setspace}
 \usepackage{esvect}
 \usepackage{graphicx}
 \usepackage[toc,page]{appendix}
 \usepackage{seqsplit}
  \usepackage{bigints}

 \usepackage[usenames]{color}
\usepackage{amssymb}
\usepackage{graphicx}
\usepackage{amscd}

\usepackage{epstopdf}% To incorporate .eps illustrations using PDFLaTeX, etc.
\usepackage[caption=false]{subfig}% Support for small, `sub' figures and tables

\usepackage[numbers,sort&compress]{natbib}% Citation support using natbib.sty
\bibpunct[, ]{[}{]}{,}{n}{,}{,}% Citation support using natbib.sty
\makeatletter% @ becomes a letter
\def\NAT@def@citea{\def\@citea{\NAT@separator}}% Suppress spaces between citations using natbib.sty
\makeatother% @ becomes a symbol again

\usepackage{authblk}

\usepackage{scalerel}[2016-12-29]
\def\stretchint#1{\vcenter{\hbox{\stretchto[440]{\displaystyle\int}{#1}}}}

\newtheorem{theorem}{Theorem}[section]

\newenvironment{proof}[1][Proof]{\begin{trivlist}
\item[\hskip \labelsep {\bfseries #1}]}{\end{trivlist}}

\newcommand{\qed}{\nobreak \ifvmode \relax \else
      \ifdim\lastskip<1.5em \hskip-\lastskip
      \hskip1.5em plus0em minus0.5em \fi \nobreak
      \vrule height0.75em width0.5em depth0.25em\fi}
\newcommand{\seqnum}[1]{\href{http://www.research.att.com/cgi-bin/access.cgi/as/~njas/sequences/eisA.cgi?Anum=#1}{\underline{#1}}}

\DeclareMathOperator{\sn}{sn}
\DeclareMathOperator{\cn}{cn}
\DeclareMathOperator{\dn}{dn}
\DeclareMathOperator{\nc}{nc}
\DeclareMathOperator{\dc}{dc}
\DeclareMathOperator{\ns}{ns}
\DeclareMathOperator{\cs}{cs}
\DeclareMathOperator{\ds}{ds}
\DeclareMathOperator{\nd}{nd}
\DeclareMathOperator{\cd}{cd}
\DeclareMathOperator{\sd}{sd}
\DeclareMathOperator{\sech}{sech}
\DeclareMathOperator{\am}{am}
\DeclareMathOperator{\gd}{gd}
\DeclareMathOperator{\arcsn}{arcsn}

\providecommand{\keywords}[1]
{
  \small	
  \textbf{\textit{Keywords---}} #1
}

\title{Exponential Riordan arrays and Jacobi elliptic functions}
\author[a]{Arnauld Mesinga Mwafise\thanks{arnauld@skku.edu}}
\author[b]{Paul Barry}
\affil[a]{Applied Algebra and Optimization Research Center, Department of Mathematics, Sungkyunkwan University, Suwon 16419, Republic of Korea}
\affil[b]{School of Science, Waterford Institute of Technology, Cork Road, Waterford, Ireland}

\begin{document}

\date{}
\maketitle

\def\x{\frac{a}{c}dP}

\def\bs{\mkern-12mu} % set amount of backspacing for lower limit of integration

\begin{abstract}
This paper establishes relationships between elliptic functions and Riordan arrays leading to new classes of Riordan arrays which here are called elliptic Riordan arrays. In particular, the case of Riordan arrays derived from Jacobi elliptic functions that are parameterized by the elliptic modulus $k$ will be treated here. Some concrete examples of such Riordan arrays are presented via a recursive formula. \end{abstract}

%TC:ignore

\section{Introduction}

Riordan arrays were originally put forward by Shapiro et al.($1991$) \cite{Shapiro} as a novel approach to representing certain classes of infinite lower triangular matrices having properties similar to the Pascal triangle in terms of their column generating functions. In addition, these Riordan arrays were found to form a group which was referred to as the Riordan group. They were so named in honour of the $20^{\text{th}}$ century combinatorialist John Riordan \cite{Riordan}. The theory of Riordan arrays was initially identified as a useful tool in solving combinatorial sum inversions \cite{Sprugnoli1, Corsani}, proof of combinatorial identities\cite{Sprugnoli2,sprugid} and in the combinatorial enumeration of lattice paths \cite{Sprugnoli1,pb1,Shapiro3}.  The two main types of Riordan arrays are the ordinary and the exponential Riordan arrays which are constructed from ordinary and exponential generating functions respectively. This paper will primarily focus on the exponential Riordan arrays. An \textit{exponential Riordan array} \cite{pbbook} is an infinite lower triangular matrix constructed from a pair of formal power series of order $0$ and $1$ respectively, represented by their  exponential generating functions given by $$g(x)=1+\sum^{\infty}_{n=1}g_n(x^n/n!) \qquad f(x)=\sum^{\infty}_{n=1}f_n(x^n/n!),$$ and with the generic element associated with the coefficients of column $k$ evaluated by $$d_{n,k}=\frac{n!}{k!}[x^n]g(x)f(x)^k.$$ The exponential Riordan array $\left(d_{n,k}\right)_{n,k\geq 0}$ is denoted by $\left[g(x),f(x)\right].$ The bivariate generating function of exponential Riordan arrays is given by \begin{equation}\label{biv2}
d_e(x,t)=\frac{1}{k!}\sum_{k=0}^{\infty}g(x)f(x)^k{t^k}=g(x)e^{tf(x)}
.\end{equation}
The case for $x=1$ in  (\ref{biv2})\; results  in the explicit formula for the row sums of the corresponding Riordan array. The most basic non-trivial exponential Riordan array is the Pascal triangle represented by $\left[e^x,x\right]$ having as its generic element $$d_{n,k}=[x^n]\frac{n!}{k!}e^{x}x^{k}\equiv{{n\choose{k}}}.$$
An alternative way of defining a Riordan array is using its recursive formula. The two main recursive rules for the formation of all the elements of a Riordan array excluding the element located on its first row are based on the $A$ and $Z$  sequence characterization of Riordan arrays \cite{deutsch,He,Rogers}. The most suitable approach of quickly determining the explicit $A$ and $Z$ generating functions of  Riordan arrays is derived from their corresponding production matrices. The production matrix $P$ of the the Riordan array $D=[g(x),f(x)]$ is defined by $$P=D^{-1}\cdot{\bar{D}}$$ where $\bar{D}$ is the matrix $D$ with its top row removed\cite{DProd,deutsch}. In particular, the Riordan array $D=\left(d_{n,k}\right)_{n,k \ge
0}$ is constructed from its  $A=\left(a_{i}\right)_{i \ge
0}$ and $Z=\left(z_{i}\right)_{i \ge
0}$ sequences as follows.
\begin{eqnarray*} (i)\qquad d_{n+1,0}&=&\sum_{i} i! z_i d_{n,i}.
\\ (ii)\qquad
d_{n+1,k}&=& a_0 d_{n,k-1}+\frac{1}{k!} \sum_{i\ge
k}i!(z_{i-k}+k a_{i-k+1})d_{n,i}. \end{eqnarray*}
The $A$ and $Z$ generating functions from the theory of production matrices are explicitly given by $$A(x)=f'(\bar{f}(x))$$ and
$$Z(x)=\frac{g'(\bar{f}(x))}{g(\bar{f}(x))}.$$ Here, $\bar{f}(x)=\text{Rev}(f)(x)$ denotes the \emph{series reversion} of $f(x)$, which is the solution $u(x)$ of the equation $f(u)=x$ that satisfies $u(0)=0$.  Furthermore, the bivariate generating function
of the matrix $P$ is given by
$$\phi_P(t,x) =
e^{tx}(Z(x)+t A(x)).$$
The Riordan array $D=\left[g(x),f(x)\right]$ can transform an infinite sequence having generating function $h(x)$ by the FTRA (Fundamental Theorem of Riordan arrays)\cite{Shapiro} to form a new infinite sequence $b(x)$ such that $$[g(x),f(x)]h(x)=g(x)h(f(x))=b(x).$$

\noindent The multiplication operation $(*)$ determines the basis for a Riordan array to form a group. The multiplication rule for two Riordan arrays $\left[g,f\right]$ and
$\left[h,l\right]$ is defined as
\begin{equation}\label{rmult}\left[g,f\right]*\left[h,l\right]=\left[g*(h\circ{f}),l\circ{f}\right].\end{equation} $[1,x]$ corresponds to  the identity element of $\left[g(x),f(x)\right]$.
 The inverse element is given by \begin{equation}\label{rinv}\left[g(x),f(x)\right]^{-1}=\left[\frac{1}{g(\bar{f}(x))},\bar{f}(x)\right].\end{equation}We note that for a power series $f(x)=\sum_{n=0}^{\infty}f_n{x^n}$ with $f(0)=0,$ we define the reversion or the compositional inverse of $f$ to be the power series $\bar{f}(x)$ such that $f(\bar{f}(x))=\bar{f}(f(x))=x.$ It can sometimes be denoted simply as $\bar{f}$ or $Rev\; f.$
The subgroups of Riordan arrays of most interest are:
Appell $\left[g(x),x\right]$, Lagrange  /associated $\left[1,f(x)\right]$
, Bell/renewal $[g(x),xg(x)]$, hitting-time $\left[\frac{xf^\prime{(x)}}{f{(x)}},f(x)\right]$, checkerboard  $[g(x),f(x)]$ where $g(x)$ is an even generating function and $f(x)$ is an odd generating function, and derivative  $[f^{\prime}(x),f(x)]$, where $f^{\prime}(x)$ denotes the first derivative of $f(x)$\cite{He}.
\\
\\
\noindent Similar to Riordan arrays, elliptic functions which are complex valued meromorphic and doubly periodic functions can be defined in terms of a reversion technique \cite{kleeh,Whit}. The origins of the main theory of elliptic functions can be traced back to the $19$th century work on integral calculus by the famous mathematician Niels H. Abel ($1802-1829$) \cite{nhabel}. His most remarkable achievement in this area was implementing the main technique of inverting elliptic integrals which led to elliptic functions. An elliptic integral can be written in the form $$\stretchint{6ex}_{\bs }^{}\!R(x,\sqrt{p(x)})dx,$$ where $R(x,w)$ is a rational function in two variables and $p(x)$ is a polynomial of degree $3$ or $4$ having no repeated roots.
During the same period when Abel put forward his work on the inversion of elliptic integrals another mathematician Carl G. Jacobi($1804-1851$) also worked in the same area. In $1829,$ Jacobi introduced the Jacobi elliptic functions denoted $\sn\;u, \cn\;u, \dn\;u.$ A key characteristic of these functions is that they satisfy the equation describing quartic elliptic curves given by
\[(y^2)^{\prime}=(1-x^2)(1-k^2{x^2}).\]
\noindent The Jacobi elliptic functions may be explicitly defined by first defining $\sn\;u\equiv{\sn(x,k)}$ such that
 \[\sn(x,k)=Rev\left(\stretchint{9ex}_{\bs 0}^x\!\frac{dt}{\sqrt{(1-t^2)(1-k^2{t^2})}}\right).\]
 We can also set $m=k^2.$ The parameter $k$ where $-1\leq{k}\leq{1}$  is called the \textit{modulus} of the elliptic integral. The \textit{complementary modulus}  is $k^{'}=\sqrt{1-k^2}$.
 The first six terms of the coefficients of the Taylor series expansion of $\arcsn(x,k) $ are
$$\left\{0,1,0,k^2+1,0,3 \left(3 k^4+2
   k^2+3\right),0\right\}.$$
 We note that we can find the power series coefficients of $\sn(x,k)$ in the first column of the inverse Riordan array

$$\left[\frac{\arcsn(x,k)}{x},  \arcsn(x,k)\right]^{-1}.$$

 \noindent An alternative method for defining $\sn(x,k)$ is to start with the function \[F(\phi,k)=\stretchint{9ex}_{\bs 0}^{\phi}\!\frac{1}{\sqrt{1-k^2{\sin(\theta)}}}d\theta.\]

 \noindent The Legendre form of the elliptic integral which is an alternative definition of elliptic $\sn(x,m)$ is given by \begin{equation}\label{ei2}arc\text{sn}(x,m)\equiv{u(\varphi,m)}=\stretchint{9ex}_{\bs 0}^{\varphi}\!\frac{d\theta}{\sqrt{1-m\sin^2{\theta}}}.\end{equation} The equivalence of these two approaches to the definition of $\text{sn}(x,m)$ follows from a change of variables which can be determined using the substitution $x=\sin{\theta}$ such that $dx=\cos(x) d\theta=\sqrt{1-x^2}d\theta$. We then revert the function $F$ to get the amplitude function $$\text{am}(u,k)=F^{-1}(u,k)=\phi.$$ Finally we define $$\sn(u,k)=\sin(\text{am}(u,k))=\sin(\phi).$$

\noindent The elliptic integral \;(\ref{ei2})\; is known as the incomplete elliptic integral. On the other hand  the complete elliptic integrals are given by

\[K(m):=u\left(\frac{\pi}{2},m\right)=\stretchint{9ex}_{\bs 0}^{1}\!\frac{dx}{\sqrt{(1-x^2)(1-mx^2)}}\]

\[u\left(\frac{\pi}{2},m\right)=\stretchint{9ex}_{\bs 0}^{\frac{\pi}{2}}\!\frac{d\theta}{\sqrt{1-m\sin^2{\theta}}}\]  with its complete complementary elliptic integral $K^{'}(m)$ defined in a similar manner but with $m$ in $K(m)$ replaced with $1-m$ such that  \[m+m^{\prime}=1\;\text{and}\; m=(k^{'})^2\;\text{and}\;\mid{m}\mid\le{1}\;\text{and}\;-K<{x}<K. \]
The basic elliptic functions satisfy the equation \[\sn^2(x;m) +\cn^2(x;m)=1.\] Thus, \[\sn(x;m)=\sqrt{1-\cn^2(x;m)}\;\;\text{and}\;
\dn(x;m)=\frac{d\varphi}{dx}=\sqrt{1-m\sn^2(x;m)}.\]
All three functions $\sn(x,k),\cn(x,k)$ and $\dn(x,k)$ are doubly periodic:
\begin{eqnarray*}
\sn(x+4K,k)&=&\sn(x,k),\;\text{where}\;\sn(K,k)=1.\\
\cn(x+4K,k)&=&\cn(x,k).\\
\dn(x+4K,k)&=&\dn(x,k).\\
\dn(x+4L,k)&=&\dn(x,k).
\end{eqnarray*}
\noindent In addition, the three basic forms of the Jacobi elliptic function determine the other 9 forms of Jacobi elliptic functions $\left(\text{sc},\nc,\dc,\ns,\cs,\ds,\nd,\cd,\sd\right)$  such that the definitions are quotients of any of these three. For example $\sd(x;m)=\frac{\sn(x;m)}{\dn(x;m)}.$
The Jacobi elliptic function are considered a generalization of the trigonometric functions from which the basic properties are defined. The basic properties of the Jacobi elliptic functions are:
\begin{itemize}
\item $\sn(0;k)=0$\qquad{$\sn(K;k)=1$}
\item $\cn(0;k)=1$\qquad{$\cn(K;k)=0$}
\item $\dn(0;k)=1$\qquad{$\dn(K;k)=k^{'}$}
\end{itemize}
In the limit we have,
\begin{eqnarray*}
\lim_{m\to{0}}\sn(x,m)&=&\sin(x)\\
\lim_{m\to{0}}\cn(x,m)&=&\cos(x)\\
\lim_{m\to{0}}\dn(x,m)&=&1\\
\lim_{m\to{1}}\sn(x,m)&=&\tanh(x)\\
\lim_{m\to{1}}\cn(x,m)&=&\sech(x)\\
\lim_{m\to{1}}\dn(x,m)&=&\sech(x).
\end{eqnarray*}

\noindent If $x=\sn(u;m)\;y=\cn(u;m)\;z=\dn(u;m)$ we get the differential system :
\[\dot{x}=yz\quad \dot{y}=-zx \quad \dot{z}=-k^2{x}y\]
satisfying the initial conditions \[\sn(0,m)=x(0)=0,\quad \cn(0,m)=y(0)=1\quad \dn(0,m)=z(0)=1.\]
\section{Elliptic functions derived from the $A$ and $Z$ generating functions of Riordan arrays}
For the inverse exponential Riordan array $M^{-1}=[g(x), f(x)]^{-1}$, we have
$$A_{M^{-1}}=\frac{1}{f'(x)},$$ and
$$Z_{M^{-1}}=-\frac{1}{f'(x)} \frac{g'(x)}{g(x)}.$$

 In particular, we can also express the array $[g, f]$ and $[g, f]^{-1}$ in terms of $\text{A}=\text{A}_M$ and $\text{Z}=\text{Z}_M$ as follows:

$$\left[\frac{1}{g(\bar{f}(x))},\bar{f}(x)\right]= \left[e^{{-}\bigintsss_{0}^x\!\frac{Z(t)}{A(t)}dt},\stretchint{6ex}_{\bs 0}^x\frac{dt}{\text{A}(t)}\right]$$
and
$$\left[\frac{1}{g(\bar{f}(x))},\bar{f}(x)\right]^{-1}=[g(x), f(x)]=\left[e^{\stretchint{4ex}_{\bs 0}^{\text{Rev}\left(\stretchint{5ex}_{\bs 0}^x\frac{dt}{\text{A}(t)}\right)} \frac{Z(t)}{\text{A}(t)}\,dt}, \text{Rev}\left(\stretchint{7ex}_{\bs 0}^x \!\frac{dt}{\text{A}(t)}\right)\right].$$

 Now we recall that an integral is called an elliptic integral if it is of the form
 $$\stretchint{7ex}_{\bs }^{}R(x, \sqrt{P(x)})\,dx,$$ where $P(x)$ is a polynomial in $x$ of degree three or four and $R$ is a rational function of its arguments.

 Thus if $$\frac{1}{A(t)}=R(t, \sqrt{P(t)}),$$ then the above exponential Riordan array satisfies the required form associated with an elliptic integral.

 As we have seen, elliptic functions are defined as the inverses of elliptic integrals. Thus the expression
$$\text{Rev}\left(\stretchint{7ex}_{\bs 0}^x\!\frac{dt}{A(t)}\right)$$ in the defining relation for $[g,f]$ is an elliptic function when
$\frac{1}{A(t)}=R(t, \sqrt{P(t)})$.

 In this case we will have
$$A(t) = \frac{1}{R(t, \sqrt{P(t)})} = \tilde{R}(t, \sqrt{P(t)}),$$ where
$$\tilde{R}=\frac{1}{R}$$ will also be a rational function.

 It is therefore natural to call an exponential Riordan array $M=[g, f]$ an \emph{elliptic Riordan array} if
$$ A_M(t)= R(t, \sqrt{P(t)}),$$ where $R$ is a rational function and $P(t)$ is a polynomial of degree three or four.

 Lagrange showed that any elliptic integral can be written in terms of the following three fundamental or normal elliptic integrals.

$$F(x,k) = \stretchint{8ex}_{\bs 0}^x \frac{dt}{\sqrt{(1-t^2)(1-k^2 t^2)}},$$

$$E(x,k) = \stretchint{8ex}_{\bs 0}^x \sqrt{\frac{1-k^2 t^2}{1-t^2}} \,dt,$$

$$\Pi(x, \alpha^2, k) = \stretchint{8ex}_{\bs 0}^x \frac{dt}{(1-\alpha^2 t) \sqrt{(1-t^2)(1-k^2t^2)}}.$$

 These integrals are called elliptic integrals of the first, second and third kind, respectively.

 As an example of elliptic functions, the Jacobi elliptic functions may be defined by first defining $\sn(x,k)$ as follows,

$$ \sn(x,k)= \text{Rev} \left(\stretchint{8ex}_{\bs 0}^x\frac{dt}{\sqrt{(1-t^2)(1-k^2 t^2)}}\right),$$

involving the elliptic integral of the first kind, and then we define

$$ \cn(x,k)= \sqrt{1- \sn(x,k)^2}\quad\text{and}\quad \dn(x,k)=\sqrt{1-k^2 \sn(x,k)^2}.$$

 Thus if $$\frac{1}{A(t)}= \frac{1}{\sqrt{(1-t^2)(1-k^2 t^2)}},$$ or
$$ A(t)= \sqrt{(1-t^2)(1-k^2 t^2)},$$ then we obtain an exponential Riordan array with
$$ f(x) = \sn(x,k).$$
\section{Jacobi Riordan arrays}
We shall refer to the first type of elliptic Riordan arrays which will be derived from Jacobi elliptic functions as Jacobi Riordan arrays. We can generate Jacobi Riordan arrays from either the formal power series of Jacobi elliptic functions or by the $A$ and $Z$ generating functions.
In general an exponential Riordan array $\left[g(x),f(x)\right]$ can be expressed in terms of its $A$ and $Z$ generating functions such that \[[g(x), f(x)]=\left[e^{\stretchint{4ex}_{\bs 0}^{\text{Rev}\left(\stretchint{4ex}_{\bs 0}^x\frac{dt}{\text{A}(t)}\right)} \frac{Z(t)}{\text{A}(t)}\,dt}, \text{Rev}\left(\stretchint{7ex}_{\bs 0}^x\frac{dt}{\text{A}(t)}\right)\right].\]
\subsection{Example:} 
Suppose that $\;A(t)=\sqrt{(1-t^2)(1-k^2{t^2})}.$
\\
\\
\noindent Then Rev$\left(\stretchint{6ex}_{\bs 0}^x \frac{dt}{\sqrt{(1-t^2)(1-k^2{t^2})}}\right)=\sn(x).$
\\
\\
\noindent So $f(x)=\sn(x).$
\\
\\
\noindent It follows that $g(x)=e^{\stretchint{4ex}_{\bs 0}^{\sn(x)}\frac{Z(t)}{A(t)}dt}.$
\\
\\
For instance, if $Z(t)=\sqrt{1-k^2{t^2}}$ then we have
  \begin{eqnarray*}
  g(x)&=&e^{\stretchint{4ex}_{\bs 0}^{\sn(x)}\frac{1}{\sqrt{1-t^2}}dt}\\
  &=&e^{\sin^{-1}(\sn(x))}
  .\end{eqnarray*}
Therefore, we have that $\left[g(x),f(x)\right]=\left[e^{\sin^{-1}(\sn(x))},\sn(x)\right].$
\\
\\
The coefficient array of $\left[e^{\sin^{-1}(\sn(x))},\sn(x)\right]$ is given by
$$\includegraphics[scale=0.7]{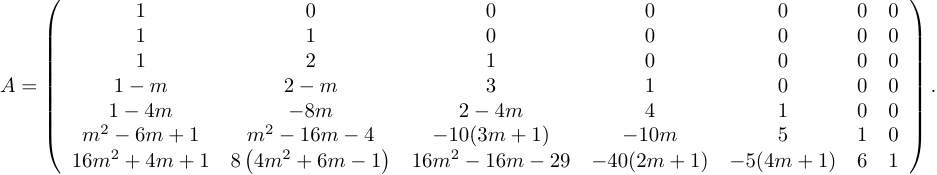}$$
\textbf{Remark 3.1.1:} The row sums of $A$ at $m=0$  corresponds to  \seqnum{A009282} in  OEIS \cite{sloane} and it also has the exponential generating function $e^{(x+\sin(x))}.$

\noindent Similarly, if $\text{Z}(t)=\sqrt{1-{t^2}}$ then we get
$$g(x)=e^{\frac{1}{k}\sin^{-1}(k \text{sn}(x))}.$$
Thus,
$$\left[g(x),f(x)\right]=\left[e^{\frac{1}{k}\sin^{-1}(k \text{sn}(x))},\text{sn}(x)\right].$$

\subsection{Example:}
Let us consider the exponential Riordan array $M=[ \cn(x), \sn(x)]$ where we suppress the parameter $k$.

We then have $$\bar{f}(x)=sn^{-1}(x).$$

We also have $g'(x)=\cn'(x)=-\sn(x) \dn(x)$, so that we obtain

\begin{eqnarray*}Z(x)&=& Z_M(x)=\frac{g'(\bar{f}(x))}{g(\bar{f}(x)}\\
&=& \frac{-\sn\left(\sn^{-1}(x)\right) \dn\left(sn^{-1}(x)\right)}{\cn\left(\sn^{-1}(x)\right)}\\
&=& \frac{- x \sqrt{1-k^2 x^2}}{\sqrt{1-x^2}}.\end{eqnarray*}

Thus for the exponential Riordan array $[\cn(x,k), \sn(x,k)]$ we have
$$ A(x)= \sqrt{(1-x^2)(1-k^2 x^2)}, \quad Z(x)=  \frac{- x \sqrt{1-k^2 x^2}}{\sqrt{1-x^2}}.$$

This means in particular that the bivariate generating function of the production matrix of $[ cn(x), sn(x)]$ is given by

$$ e^{xy}\left(\frac{- x \sqrt{1-k^2 x^2}}{\sqrt{1-x^2}}+y \sqrt{(1-x^2)(1-k^2 x^2)}\right).$$

In this case, we have

$$\frac{Z(t)}{A(t)} = \frac{- t \sqrt{1-k^2 t^2}}{\sqrt{1-t^2}} \frac{1}{\sqrt{(1-x^2)(1-k^2 t^2)}}$$ or
$$\frac{Z(t)}{A(t)} = \frac{-t}{1-t^2}.$$

Thus

\begin{eqnarray*}e^{\stretchint{4ex}_{\bs 0}^{\text{Rev}\left(\stretchint{4ex}_{\bs 0}^x\frac{dt}{A(t)}\right)} \frac{Z(t)}{A(t)}\,dt}&=&
e^{\stretchint{4ex}_{\bs 0}^{\sn(x,k)}\!\frac{-t}{1-t^2} \,dt}\\
&=& \sqrt{1- \sn(x,k)^2}\\
&=& \cn(x,k),\end{eqnarray*}
as expected.

To calculate the inverse array $[\cn(x), \sn(x)]$ we have

$$\frac{1}{g(\bar{f}(x))}= \frac{1}{\cn(\sn^{-1}(x))}=\frac{1}{\sqrt{1-x^2}},$$ and so
$$[\cn(x), \sn(x)]^{-1} = \left[\frac{1}{\sqrt{1-x^2}}, \stretchint{8ex}_{\bs 0}^x \frac{dt}{\sqrt{(1-k^2 t^2)(1-t^2)}}\right].$$

\subsection{Example:}

We next look at the exponential Riordan array

$$\left[\frac{\cn(x)}{1+\sn(x)}, \sn(x)\right].$$
We once again have
$$A(x)= \sqrt{(1-x^2)(1-k^2 x^2)}.$$
Now $g(x)=\frac{\cn(x)}{1+\sn(x)}$, and we find that
$$g'(x)=\frac{-(1+\sn(x))\sn(x)\dn(x)-\cn(x)^2 \dn(x)}{(1+\sn(x))^2}.$$

We then get
\begin{eqnarray*}
\frac{g'(\sn^{-1}(x)}{g(\sn^{-1}(x))}&=&\frac{-(1+x)x\sqrt{1-k^2 x^2}-(1-x^2)\sqrt{1-k^2 x^2}}{(1+x)^2} \frac{1+x}{\sqrt{1-x^2}}\\
&=& -\frac{\sqrt{1-k^2x^2}}{\sqrt{1-x^2}}.\end{eqnarray*}
Thus the bivariate generating function for the production matrix of $\left[\frac{\cn(x)}{1+\sn(x)}, \sn(x)\right]$ is given by

$$e^{xy}\left(-\sqrt{\frac{1-k^2x^2}{1-x^2}}+y \sqrt{(1-x^2)(1-k^2 x^2)}\right).$$

We note that for $k^2=1$, we get the exponential Riordan array
$$\left[\frac{\sech(x)}{1+\tanh(x)}, \tanh(x)\right].$$

The inverse matrix is calculated as follows.
\begin{eqnarray*}
\left[\frac{\cn(x)}{1+\sn(x)}, \sn(x)\right]^{-1}&=& \left[\frac{1}{\frac{\cn(\sn^{-1}(x)}{1+\sn(\sn^{-1}(x))}}, sn^{-1}(x)\right]\\
&=& \left[\frac{1+x}{\sqrt{1-x^2}}, \stretchint{8ex}_{\bs 0}^x \frac{dt}{\sqrt{(1-k^2t^2)(1-t^2)}}\right].\end{eqnarray*}

\subsection{Example:}
Our next example is the exponential Riordan array
$$\left[\frac{\cn(x)}{1+\sn(x)}, \frac{\sn(x)}{1+\sn(x)}\right].$$

We do not immediately know what the inverse function of $\frac{\sn(x)}{1+\sn(x)}$ is, so we use the theory of Riordan arrays to continue the analysis.

Thus we note that

$$\left[\frac{\cn(x)}{1+\sn(x)}, \frac{\sn(x)}{1+\sn(x)}\right]=[\cn(x), \sn(x)] \left[\frac{1}{1+x}, \frac{x}{1+x}\right],$$ where the second Riordan array in the product is related to the Laguerre polynomials.

Taking inverses, we obtain

$$\left[\frac{\cn(x)}{1+\sn(x)}, \frac{\sn(x)}{1+\sn(x)}\right]^{-1} = \left[\frac{1}{1+x}, \frac{x}{1+x}\right]^{-1} [\cn(x), \sn(x)]^{-1},$$

or

$$\left[\frac{\cn(x)}{1+\sn(x)}, \frac{\sn(x)}{1+\sn(x)}\right]^{-1}=\left[\frac{1}{1-x}, \frac{x}{1-x}\right]\left[\frac{1}{\sqrt{1-x^2}},\stretchint{7ex}_{\bs 0}^x \frac{dt}{(1-t^2)(1-k^2t^2)}\right].$$

This gives us

$$\left[\frac{\cn(x)}{1+\sn(x)}, \frac{\sn(x)}{1+\sn(x)}\right]^{-1}=\left[\frac{1}{\sqrt{1-2x}}, \stretchint{8ex}_{\bs 0}^{\frac{x}{1-x}} \!\frac{dt}{(1-t^2)(1-k^2t^2)}\right].$$

Thus in particular, we have that

$$\frac{\sn(x)}{1+\sn(x)}=\text{Rev} \stretchint{7ex}_{\bs 0}^{\frac{x}{1-x}} \!\frac{dt}{(1-t^2)(1-k^2t^2)},$$

or by the change of variable $y=\frac{t}{1+t}$, we get

$$ \frac{\sn(x)}{1+\sn(x)}=\text{Rev}\stretchint{8ex}_{\bs 0}^x\frac{dy}{\sqrt{(1-2y)(1-2y-(k^2-1)y^2)}}.$$

We can generalize this to the following.

\begin{theorem}
Given two exponential Riordan arrays, if any one of these Riordan arrays is an elliptic Riordan array, then their product is also an elliptic Riordan array.
\end{theorem}
\begin{proof}
We assume given two exponential Riordan arrays, $M=[ g, f]$, and $N=[u, v]$. We assume that
$M$ is an elliptic Riordan array, with
$$A_M(t)= R(t, \sqrt{P(t)}).$$

We consider the product
$$M\cdot N = [g, f] \cdot [u, v] = [g u(f), v(f)].$$

Knowing that $$f(x)=\text{Rev} \stretchint{7ex}_{\bs 0}^x \frac{dt}{A(t)},$$ we wish to find an ``elliptic'' characterisation of
$v(f)$.

For this, we look at the inverse

$$ (M \cdot N)^{-1} = N^{-1}\cdot M^{-1} = \left[\frac{1}{u(\bar{v})}, \bar{v}\right]\cdot \left[\frac{1}{g(\bar{f})}, \bar{f}\right].$$

Now $$\bar{f}(x)=\stretchint{7ex}_{\bs 0}^x \frac{dt}{A(t)},$$ so we obtain

$$ (M \cdot N)^{-1}= \left[\frac{1}{u(\bar{v})} \frac{1}{g(\bar{f}(\bar{v}))}, \stretchint{7ex}_{\bs 0}^{\bar{v}(x)} \!\frac{dt}{A(t)}\right].$$

Thus we have that

$$ v(f) = \text{Rev}\stretchint{7ex}_{\bs 0}^{\bar{v}(x)} \!\frac{dt}{A(t)}.$$

To put this in an ``elliptic'' form, we use the change of variable

$$ y = v(t) \Longrightarrow  t=\bar{v}(y). $$

in the integral.

This gives us

$$ \frac{dy}{dt} = v'(t) \Longrightarrow dt = \frac{dy}{v'(t)} = \frac{dy}{v'(\bar{v}(y))} = \bar{v}'(y)dy.$$

When $t = \bar{v}(x)$, we have $y=v(t)=v(\bar{v}(x))=x$, and so we have

$$\stretchint{7ex}_{\bs 0}^{\bar{v}(x)} \!\frac{dt}{A(t)} = \stretchint{7ex}_{\bs 0}^x \frac{\bar{v}'(y)dy}{A(\bar{v}(y))}.$$

Thus we have

$$ v(f) = \text{Rev}\stretchint{8ex}_{\bs 0}^x \frac{\bar{v}'(y)dy}{A(\bar{v}(y))}
= \text{Rev}\stretchint{8ex}_{\bs 0}^x \frac{dy}{v'(\bar{v}(y)) A(\bar{v}(y))}.$$

\end{proof}

\subsection{Example:} In this example, we seek to write the elliptic function
$$\frac{\sn(x)(1+\sn(x))}{1-\sn(x)}$$ as the reversion of an integral whose limits are $0$ to $x$.
For this, we consider the elliptic Riordan array $[\cn(x), \sn(x)]$ and the transformation given by the exponential Riordan array $\left[\frac{1}{1-x}, \frac{x(1+x)}{1-x}\right]$.

Thus we have the product

$$[\cn(x), \sn(x)] \cdot \left[\frac{1}{1-x}, \frac{x(1+x)}{1-x}\right]= \left[\frac{cn(x)}{1-sn(x)}, \frac{\sn(x)(1+\sn(x))}{1-\sn(x)}\right].$$

Here, $$\left[\frac{1}{1-x}, \frac{x(1+x)}{1-x}\right]^{-1}=\left[1-\frac{\sqrt{1+6x+x^2}-x-1}{2}, \frac{\sqrt{1+6x+x^2}-x-1}{2}\right].$$

In particular, $$\bar{v}(x)=\frac{\sqrt{1+6x+x^2}-x-1}{2},$$ and
$$\bar{v}'(x)=\frac{3+x-\sqrt{1+6x+x^2}}{2 \sqrt{1+6x+x^2}}.$$

Also, $$sn(x) = \text{Rev}\stretchint{8ex}_{\bs 0}^x \frac{dt}{(1-k^2 t^2)(1-t^2)}.$$

Thus we have
$$ \frac{\sn(x)(1+\sn(x))}{1-\sn(x)}= \text{Rev}\stretchint{18ex}_{\bs 0}^x \frac{\frac{3+y-\sqrt{1+6y+y^2}}{2 \sqrt{1+6y+y^2}}dy}{\sqrt{\left(1-k^2 \left(\frac{\sqrt{1+6y+y^2}-y-1}{2}\right)^2\right)\left(1-\left(\frac{\sqrt{1+6y+y^2}-y-1}{2}\right)^2\right)}}.$$

In the next examples that follow we examine some Jacobi Riordan arrays based on their coefficient and production matrices representations.

\subsection{Example: $\left[\dn(x,m),\sn(x,m)\right]$ } The coefficient matrix of $\left[\dn(x,m),\sn(x,m)\right]$ begins
$$\includegraphics[scale=0.5]{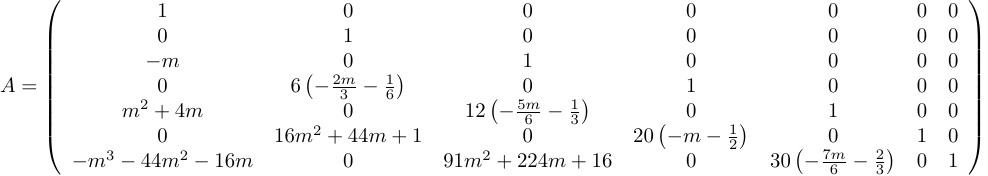}$$ which is equivalent to $$\includegraphics[scale=0.5]{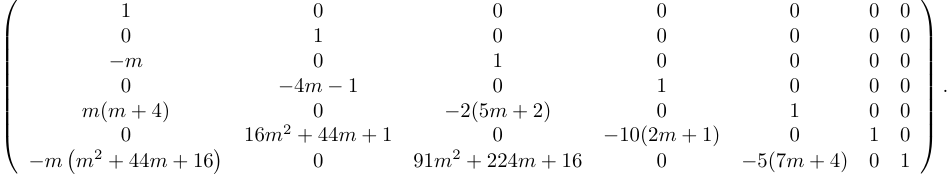}$$
\noindent \textbf{Remark 3.6.1:} 
\begin{itemize}
\item The numbers $1,6,12,20,30,...$ positioned along the $n+2,n$ diagonal of  the matrix $A$ corresponds to OEIS \seqnum{A180291.}
\item The row sums of $A$ for $m=0$ form the sequence $\left(1, 1, 1, 0, -3, -8, -3,... \right)$ which correspond to OEIS \seqnum{A002017}  with exponential generating function  \; $e^{\sin(x)}.$
\item The row sums of $A$ for $m=1$ form the sequence $\left(1, 1, 0, -4, -8, 32, 216,... \right)$ which  corresponds to OEIS \seqnum{A009265} with exponential generating function \;$\frac{e^{\tanh(x)}}{cosh(x)}.$
\end{itemize}
The production matrix  of $A$ in terms of $m$ is given by:
$$\includegraphics[scale=0.5]{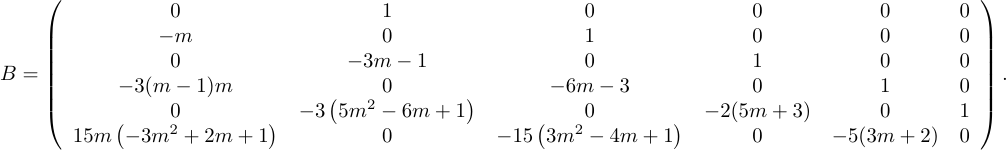}$$
If $m=-1,0,1$ then $\left[\text{dn}(x,m),\text{sn}(x,m)\right]$
 produces the Riordan arrays $$C=\{\left[\text{dn}(x,-1),\text{sn}(x,-1)\right],\left[1,\sin(x)\right],\left[\text{sech}(x),\tanh(x)\right]\}\;\text{respectively}.$$
\noindent The production matrices of $B$ associated to the Riordan arrays in $C$ for $m=-1,0,1$ are as follows:
$$\includegraphics[scale=0.5]{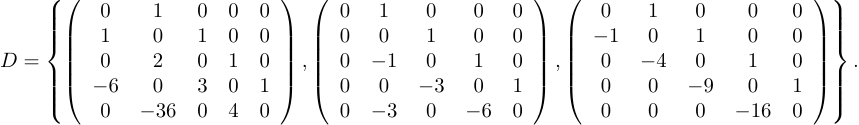}$$
\noindent\textbf{Remark 3.6.2:} The tri-diagonal matrix for $m=1$ in $D$ indicates that the inverse of the Riordan array $\left[\text{sech}(x),\tanh(x)\right]$ forms the coefficient array of a family of orthogonal polynomials \cite{Bhankel}. This is OEIS \seqnum{A060524} the number of degree $n$-permutations with $k$ odd cycles.
$\left[\frac{1}{\sqrt{1-x^2}},\tanh ^{-1}(x)\right]$ is the coefficient array of a family of orthogonal polynomials. The three term recurrence relation for the family of orthogonal polynomials is given by \[P_{n+1}(x)=xP_n(x)+n^2{P_{n-1}}(x),\;\forall{n\geq{1}}\] with $P_0(x)=1,\;P_1(x)=x.$ $$\text{In particular, let}\;Q_n(x)=\frac{P_n(ix)}{i^n},\;(i^2=-1)$$ we get $$Q_{n+1}(x)=xQ_{n}(x)-n^2{Q_{n-1}}(x),\;\forall{n\geq{1}}.$$
\noindent The  $A$ and $Z$ sequences associated to the production matrices of $B$ and $D$ are listed below:
\begin{eqnarray*}
A(x,m)&=&\text{cn}\left(\left.\text{sn}^{-1}(x,m)\right|m\right)\text{dn}\left(\left.\text{sn}^{-1}(x,m)\right|m\right)\\
&=&(\sqrt{1-x^2})(\sqrt{1-mx^2})\\
&=&\sqrt{(1-x^2)(1-mx^2)}.\\
 A(x,1)&=&1-x^2.\\
 A(x,0)&=&\sqrt{1-x^2}.\\
 Z(x,m)&=&\frac{-m\text{sn}(\bar{\text{sn}}(x,m),m)\text{cn}(\bar{\text{sn}}(x,m),m)}{\text{dn}(\bar{\text{sn}}(x,m),m)}\\
 &=&\frac{-mx\sqrt{1-x^2}}{\sqrt{1-mx^2}}.\\
Z(x,1)&=&-x.
\end{eqnarray*}
\subsection{Example: $\left[\frac{d}{dz}\sn(x,m),\sn(x,m)\right]$ }
The  coefficient array of $\left[\frac{d}{dx}\sn(x,m),\sn(x,m)\right]$ where $\frac{d}{dx}\sn(x,m)=\cn(x|m) \dn(x|m)$ is given by:
$$\includegraphics[scale=0.7]{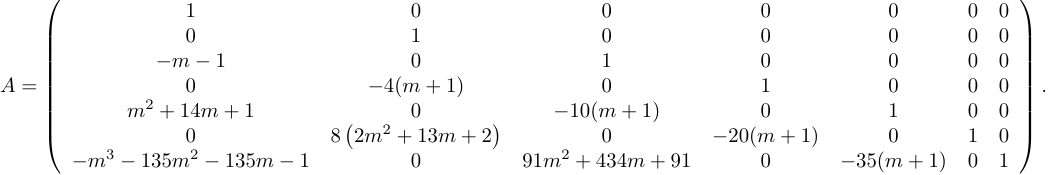}$$
\noindent\textbf{Remark 3.7.1:}
\begin{itemize}
\item The coefficient matrix of $\left[\frac{d}{dx}\sn(x,m),\sn(x,m)\right]$ forms a \textit{palindromic} Riordan array.

\item The row sums of $A$ for $m=0$ form the sequence $\left(1, 1, 0, -3, -8, -3, 56,...\right)$  which has the exponential generating function \;$\cos(x)e^{\sin(x)}.$

\item The row sums of $A$ for $m=1$ form the sequence $\left(1, 1, -1, -7, -3, 97, 275,... \right)$  which has the exponential generating function \;$\text{sech}^2(x)e^{\tanh(x)}.$
\item The non-zero entries of the first column of $A$ corresponds to the matrix $$\includegraphics[scale=0.7]{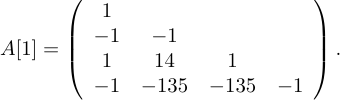}$$
By multiplying the matrix $A[1]$ by $-1^{m}$ if $n\equiv{m}(\textrm{mod{\;2}})$ where $n$ is the column number s.t $n=0,1,2,...$ we get the matrix $$\includegraphics[scale=0.7]{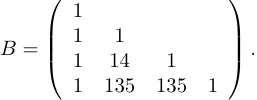}$$
The matrix $B$ corresponds to OEIS \seqnum{A060628.}
\end{itemize}
\noindent The production matrix of $A$ in terms of $m$ is given by
$$\includegraphics[scale=0.7]{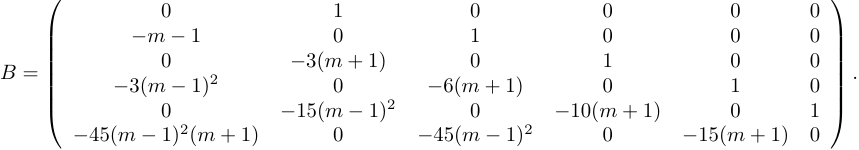}$$
\\
\noindent If $m=-1,0,1$ then  $\left[\frac{d}{dx}\sn(x,m),\sn(x,m)\right]$ produces the Riordan arrays
$$C=\{\left[\cn(x,-1) \dn(x,-1),\sn(x,-1)\right],\left[\cos (x),\sin (x)\right],\left[\sech^2(x),\tanh(x)\right]\}\;\text{respectively.}$$
\\
\\
\noindent The production matrices of  $B$  in terms of $m=-1,0,1$ are as follows:

$$\includegraphics[scale=0.7]{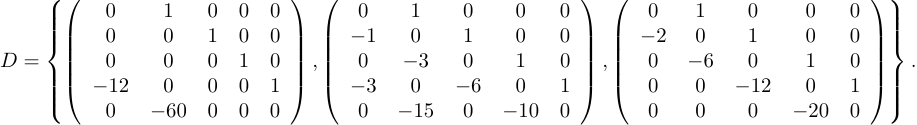}$$
\textbf{Remark 3.7.2:} The tridiagonal production matrix for $m=1$ in $D$ which is associated to the Riordan array $\left[\sech^2(x),\tanh(x)\right]$ in $C$ forms an orthogonal polynomial sequence for $\left[\sech^2(x),\tanh(x)\right]^{-1}.$
 Furthermore, \[\left[\sech^2(x),\tanh(x)\right]^{-1}=\left[\frac{1}{1-x ^2},\tanh ^{-1}(x )\right]\]
 represents the coefficient matrix of the family of orthogonal polynomials. The three term recurrence relation for these polynomials is given by
 \[P_{n+1}(x)=x{P}_n(x)+n(n+1)P_{n-1}(x)\] with $P_0(x)=1,\;P_1(x)=x\;\;\text{s.t.}\;-1<x <1.$
\\
\\
\noindent The  $A$ and $Z$ sequences corresponding to the production matrix of  $B$ and $D$ are listed as follows:
\begin{eqnarray*}
A(x,m)&=&\text{cn}\left(\left.\text{sn}^{-1}(x|m)\right|m\right)
 \text{dn}\left(\left.\text{sn}^{-1}(x|m)\right|m\right)\\
 &=&\sqrt{(1-sn^2(\bar{sn}(x,m),m)(1-m^2{sn^2(\bar{sn}(x,m),m)})}\\
 &=&\sqrt{(1-x^2)(1-m^2{x^2})}\\
A(x,1)&=&1-x^2\\
 A(x,0)&=&\sqrt{1-x^2}\\
 Z(x,m)&=&\frac{x \left(m \left(2 x^2-1\right)-1\right)
   \text{cn}\left(\left.\text{sn}^{-1}(x|m)\right|m\right)
   \text{dn}\left(\left.\text{sn}^{-1}(x|m)\right|m\right)}{\left(x^2-1\right) \left(m x^2-1\right)}\\
   &=&\frac{x \left(m \left(2 x^2-1\right)-1\right)\sqrt{(1-sn^2(\bar{sn}(x,m),m)(1-m^2{sn^2(\bar{sn}(x,m),m)})}}{\left(x^2-1\right) \left(m x^2-1\right)}\\
   &=&\frac{x \left(m \left(2 x^2-1\right)-1\right)\sqrt{(1-x^2)(1-m^2{x^2})}}{\left(x^2-1\right) \left(m x^2-1\right)}\\
   Z(x,1)&=&\frac{-x\left(2
   x^2-2\right)}{\left(x^2-1\right)}\\
   &=&-2x\\
   Z(x,0)&=&-\frac{x}{\sqrt{1-x^2}}.
  \end{eqnarray*}

  \subsection{$\left[\frac{d}{dx}\text{am}(x,m),\text{am}(x,m)\right]$} The coefficient array $\left[\frac{d}{dz}\text{am}(x,m),\text{am}(x,m)\right]$ where $\frac{d}{dx}\text{am}(x,m)=\text{dn}(x,m)$
is given by

$$\includegraphics[scale=0.7]{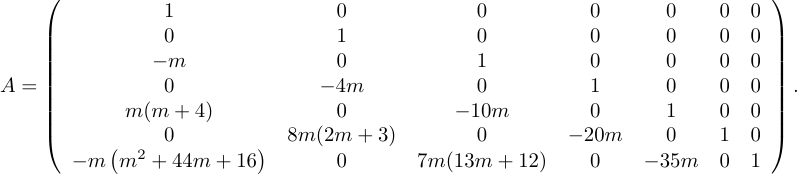}$$
\textbf{Remark 3.8.1:}
\begin{itemize}
\item For $m=0$ we have $A=I=\left[1,x\right].$
\item The row sums of $A$ for $m=1$ form the sequence $\left(1, 1, 0, -3, -4, 21, 80,... \right)$  corresponding to OEIS    \seqnum{A012123} with exponential generating function 

$e^{\sin^{-1}(\tanh(x))}=e^{\gd(x)}$ where $\gd(x)$ is the \textit{Gudermannian function} such that \[\gd(x)=\stretchint{7ex}_{\bs 0}^x \frac{1}{\cosh{t}}\,dt\qquad{-\infty<x<\infty}.\]
\item The nonzero elements of the first column of the matrix $A$ generated from the derivative of the Jacobi amplitude function  forms the coefficient matrix $$\includegraphics[scale=0.7]{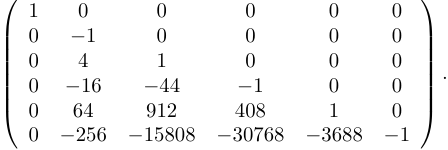}$$
\end{itemize}
\noindent The production matrix of $A$ in terms of $m$:
$$\includegraphics[scale=0.7]{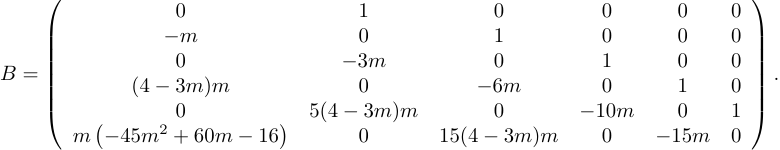}$$
\noindent If $m=1$ and $m=0$ then $\left[\frac{d}{dx}\am(x,m),\am(x,m)\right]$ produces the Riordan arrays $$C=\{\left[\frac{2 e^x}{e^{2 x}+1},2 \tan ^{-1}\left(e^x\right)-\frac{\pi }{2}\right],\left[1,x\right]\}\;\text{respectively.}$$
\noindent\textbf{Remark 3.8.2:} We note that $\frac{2 e^x}{e^{2 x}+1}$ has an ordinary generating function given by
\[\cfrac{1}{1+\cfrac{x^2}{1+\cfrac{4x^2}{1+
\cfrac{9 x^2}{1+
\cfrac{25 x^2}{1+\cdots}}}}}.\]
\noindent The production matrices from $B$ in terms of $m=-1,0,1$ are as follows:
$$\includegraphics[scale=0.7]{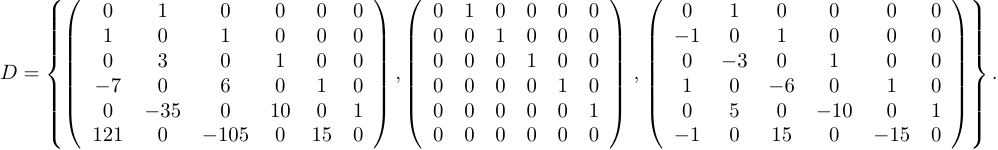}$$
 \textbf{Remark 3.8.3:} We note that the generating function of the matrix $D$ at $m=1$ is $$e^{xy}(-\sin(x)+y\cos(x)).$$
\noindent The  $A$ and $Z$ sequences corresponding to the production matrices of  $B$ and $D$ are listed as follows:
\begin{eqnarray*}
A(x,m)&=&\text{dn}(F(x,m),m)\\
  A(x,1)&=&\text{sech}(F(x,1))\\
  &=&\text{sech} (\log (\tan (x)+\sec (x)))\\
  &=&\frac{2 (\tan (x)+\sec (x))}{(\tan (x)+\sec
   (x))^2+1}\\
   &=&\cos(x)\\
  Z(x,m)&=&\frac{m \text{cn}(F(x,m),m)
   \text{sn}(F(x,m),m)}{\text{dn}(F(x,m),m)}\\
   Z(x,1)&=&-\tanh (F(x,1))\\
   &=&-\tanh(\log (\tan (x)+\sec (x)))\\
   &=&\frac{1-(\tan (x)+\sec (x))^2}{1+(\tan (x)+\sec
   (x))^2}\\
   &=&\sin(x).\\
    \end{eqnarray*}

\section{Conclusion} 

We have shown that through the definition of exponential Riordan arrays by means of their A and Z sequences, that it is natural to consider 
those exponential Riordan arrays that are defined by elliptic functions, when those elliptic functions are defined by a reversion. We have considered a number of such arrays defined by common Jacobi elliptic functions. In further work, we intend to re-cast some applications of elliptic functions to this new context, thus gaining more insight into the use of Riordan arrays in applied problems, as well as augmenting the solution context. Candidate areas for such explorations include transmission line theory, soliton theory and cosmology.

\end{document}